\documentclass[10pt]{amsart}
\pdfoutput=1

\usepackage[text={120mm,178mm},centering]{geometry}
\usepackage[usenames,dvipsnames]{xcolor}
\usepackage[colorlinks, 
            linkcolor=MidnightBlue, 
            citecolor=MidnightBlue,
            urlcolor =BlueViolet,
           ]{hyperref}

\usepackage{titlesec}

\titleformat{\section}[hang]
  {\normalfont\sc\filcenter}{\thesection.}{0.3em}{}

\titleformat{\subsection}[runin]
               {\normalfont\bfseries}
               {{\rm{\S}}\thesubsection.}{.4em}{}[.]

\makeatletter

\def\@setauthors{%
  \begingroup
  \def\thanks{\protect\thanks@warning}%
  \trivlist
  \centering\footnotesize \@topsep30\p@\relax
  \advance\@topsep by -\baselineskip
  \item\relax
  \author@andify\authors
  \def\\{\protect\linebreak}%
  {\normalsize\textsc\authors}%
  \ifx\@empty\contribs
  \else
    ,\penalty-3 \space \@setcontribs
    \@closetoccontribs
  \fi
  \endtrivlist
  \endgroup
}

\renewenvironment{abstract}{%
  \ifx\maketitle\relax
    \ClassWarning{\@classname}{Abstract should precede
      \protect\maketitle\space in AMS document classes; reported}%
  \fi
  \global\setbox\abstractbox=\vtop \bgroup
    \normalfont\Small
    \list{}{\labelwidth\z@
      \leftmargin8.5mm \rightmargin\leftmargin
      \listparindent\normalparindent \itemindent\z@
      \parsep\z@ \@plus\p@
      
    }%
    \item[\hskip\labelsep\scshape\abstractname.]%
}{%
  \endlist\egroup
  \ifx\@setabstract\relax \@setabstracta \fi
}
\def\@setabstract{\@setabstracta \global\let\@setabstract\relax}
\def\@setabstracta{%
  \ifvoid\abstractbox
  \else
    \skip@20\p@ \advance\skip@-\lastskip
    \advance\skip@-\baselineskip \vskip\skip@
    \box\abstractbox
    \prevdepth\z@ 
  \fi
}

\def\@biblabel#1{\@ifnotempty{#1}{#1.}}

  \def\@typesizes{%
    \or{5}{6}\or{6}{7}\or{8}{8}\or{8}{10}\or{9}{11}%
    \or{10}{12}
    \or{\@xipt}{13}\or{\@xiipt}{14}\or{\@xivpt}{17}%
    \or{\@xviipt}{20}\or{\@xxpt}{24}}%

\def\ps@firstpage{\ps@plain
  \def\@oddfoot{\normalfont\scriptsize \hfil\thepage\hfil
     \global\topskip\normaltopskip}%
  \let\@evenfoot\@oddfoot
  \def\@oddhead{{\scriptsize\llap{\fontfamily{ptm}\selectfont Progress} {\fontfamily{ptm}\selectfont in Mathematics, Vol.~195, \copyright\ 2001 Birkh\"auser Verlag Basel/Switzerland}}\hfil}%
  \let\@evenhead\@oddhead 
}

\def\@makefntext{\indent\hskip 3.5mm\@makefnmark}

\makeatother

\newcommand\fancybreak{\goodbreak\medskip\centerline{$\circ$\quad\raise 0.5ex\hbox to 1in{\hrulefill}\quad\raise 0.1ex\hbox to 1em{$\scriptstyle\diamondsuit$}\quad\raise 0.5ex\hbox to 1in{\hrulefill}\quad$\circ$}\medskip\goodbreak}

\newtheorem{theorem}{Theorem}[section]
\newtheorem{lemma}[theorem]{Lemma}
\newtheorem{proposition}[theorem]{Proposition}

\newcommand\CC{{\mathcal C}}
\newcommand\CBtwo{{\mathcal B}_{2}}
\newcommand\CBone{{\mathcal B}_{1}}

\newcommand\BGm{{\mathbf G_{m}}}
\newcommand\BQ{{\mathbf Q}}
\newcommand\BZ{{\mathbf Z}}

\newcommand\pid{{\mathfrak p}}

\newcommand\ra{\rightarrow}
\newcommand\kb{k^{\textup{sep}}}
\newcommand\Xh{\hat{X}}
\newcommand\Yh{\hat{Y}}
\newcommand\Qb{\bar{Q}}
\newcommand\db{\discretionary{}{}{}}

\let\hat=\widehat
\let\bar=\overline

\DeclareMathOperator\Aut{Aut}
\DeclareMathOperator\End{End}
\DeclareMathOperator\Gal{Gal}
\DeclareMathOperator\inv{inv}
\DeclareMathOperator\Ker{Ker}
\DeclareMathOperator\Nrd{Nrd}
\DeclareMathOperator\Prin{Prin}
\DeclareMathOperator\rank{rank}
\DeclareMathOperator\Res{Res}
\DeclareMathOperator\Spec{Spec}
\DeclareMathOperator\Trd{Trd}

\begin{document}

\title[Abelian varieties with no principal polarizations]
      {Isogeny classes of abelian varieties\\ with no principal polarizations\\\ }
\author{Everett W. Howe}

\keywords{Abelian variety, isogeny, polarization, Galois twist}
\subjclass[1991]{Primary 14K02; Secondary 11G10, 14K15}


\begin{abstract}
We provide a simple method of constructing isogeny classes of abelian varieties
over certain fields $k$ such that no variety in the isogeny class has a 
principal polarization. In particular, given a field~$k$, a Galois extension
$\ell$ of $k$ of odd prime degree~$p$, and an elliptic curve $E$ over $k$ that
has no complex multiplication over $k$ and that has no $k$-defined $p$-isogenies
to another elliptic curve, we construct a simple $(p-1)$-dimensional abelian
variety $X$ over $k$ such that every polarization of every abelian variety
isogenous to $X$ has degree divisible by~$p^{2}$. We note that for every odd
prime $p$ and every number field~$k$, there exist $\ell$ and $E$ as above. We
also provide a general framework for determining which finite group schemes
occur as kernels of polarizations of abelian varieties in a given isogeny class.

Our construction was inspired by a similar construction of Silverberg and 
Zarhin; their construction requires that the base field $k$ have positive 
characteristic and that there be a Galois extension of $k$ with a certain 
non-abelian Galois group.
\end{abstract}

\thispagestyle{empty}
\phantom{.}
\bigskip\bigskip

\centerline{Corrigendum to}
\bigskip
\centerline{Everett W.~Howe: Isogeny classes of abelian varieties with no
            principal polarizations,}
\centerline{pp. 203--216 in: \emph{Moduli of Abelian Varieties} 
            (C.~Faber, G.~van der Geer, F.~Oort, eds.),}
\centerline{Progress in Mathematics \textbf{195}, Birkh\"auser, Basel, 2001}

\bigskip
\fancybreak
\bigskip\bigskip\bigskip

The purpose of this note is to point out an error in my paper
\emph{Isogeny classes of abelian varieties with no principal polarizations},
which was originally put on the arXiv in February 2000. I have typeset this note
so that it closely matches the formatting of the original paper as it appeared
in the 2001 volume \emph{Moduli of Abelian Varieties}, and I have highlighted
portions of the text in red and added explanatory comments in the margin.

Note that the material in the first two sections is correct. The errors 
presented here all occur in Section~3. In particular, the proof of Theorem~3.2
is incorrect. The ``only if'' implication in the theorem holds (with the
choice of the set $S_{\CC}$ given in the proof), but the ``if'' implication does
not. 

The error in the proof comes down to the incorrect claim that if $H$ is a finite
simple group scheme over a field $k$ and if there is a nondegenerate alternating 
pairing $e$ on $H\times H$, then there is an embedding $H\to H\times H$ whose
image is isotropic with respect to $e$. This claim is true when $k$ is finite,
but false, for example, when $k=\BQ$. A counterexample can be given by taking 
$H$ to be the $3$-torsion of an elliptic curve over $\BQ$ for which the Galois
action is maximal, and taking $e$ to be the product of the Weil pairings on the
factors of~$H\times H$.

\bigskip

I am grateful to Armand Brumer for bringing this issue to my attention and for
suggesting the counterexample to the claim.

\bigskip\bigskip

\hskip 0.5\hsize --- Everett W. Howe

\hskip 0.5\hsize\phantom{---} San Diego, California

\hskip 0.5\hsize\phantom{---} 10 December 2022

\vfill\eject

\vbox{\renewcommand{\baselinestretch}{0.8}\footnotesize\noindent\color{Red}
Note: I have no record of signing a copyright transfer agreement for this work.
But for verisimilitude, I have included the copyright notice above, which
appears in the published version of the paper. --- EWH}

\vskip 2ex

\maketitle
\phantom{.}
\vspace{-5ex}

\setcounter{page}{203}

\section{Introduction}
\label{PARTintroduction}

\noindent
A natural question to ask of an isogeny class $\CC$ of abelian varieties over a
field $k$ is whether or not it contains a principally polarized variety. If $k$
is algebraically closed then $\CC$ will certainly contain a principally
polarized variety, and if $k$ is finite then every $\CC$ that satisfies some
relatively weak conditions will contain a principally polarized variety
(see~\cite{brefHoweAV}); for example, it is enough for the varieties in $\CC$ to
be simple and odd-dimensional. In this paper we show that for a large class of
fields, including all number fields and function fields over finite fields, it 
is very easy to construct isogeny classes of abelian varieties that contain no
principally polarized varieties. We also provide a framework for considering the
more general problem of determining which finite group schemes occur as the
kernels of polarizations of varieties in a given isogeny class.

Our construction of isogeny classes containing no principally polarized
varieties is very straightforward, but to describe it we must introduce some
terminology. If $\ell$ is a finite extension of a field $k$ and if $E$ is an 
elliptic curve over~$\ell$, we let $\Res_{\ell/k} E$ denote the restriction of
scalars of $E$ from $\ell$ to $k$ (see Section~1.3 of~\cite{brefWeil}). If $E$
is an elliptic curve over~$k$, we define the 
\emph{reduced restriction of scalars} of $E$ from $\ell$ to $k$ to be the kernel
of the trace map from $\Res_{\ell/k} E$ to~$E$. Let $p$ be an odd prime. We will
say that an elliptic curve $E$ over $k$ is \emph{$p$-isolated} if it has no
complex multiplication over $k$ and if it has no $p$-isogeny to another elliptic
curve over~$k$. We will say that a field $k$ is \emph{$p$-admissible} if there 
is a $p$-isolated elliptic curve over $k$ and if there is a Galois extension of
$k$ of degree~$p$. One can show, for example, that every number field is 
$p$-admissible; that every function field over a finite field of characteristic
not $p$ is $p$-admissible; and that no finite field or algebraically-closed
field is $p$-admissible.

\begin{theorem}
\label{TTmain}
Let $p$ be an odd prime number and let $k$ be a $p$-admissible field. Let $E$ be
a $p$-isolated elliptic curve over~$k$, let $\ell$ be a degree-$p$ Galois 
extension of~$k$, and let $X$ be the reduced restriction of scalars of $E$ from
$\ell$ to~$k$. Then $X$ is simple, and every polarization of every abelian
variety isogenous to $X$ has degree divisible by~$p^{2}$.
\end{theorem}

We provide an elementary proof of this theorem in Part~\ref{PARTnoPPAVs} of this
paper. In Section~\ref{SSgeneral} we prove some basic results about 
polarizations and endomorphism rings of Galois twists of abelian varieties; we 
apply these general results in Section~\ref{SSrestriction} to prove several
results about reduced restrictions of scalars of $p$-isolated elliptic curves. 
We use the results of Section~\ref{SSrestriction} to prove Theorem~\ref{TTmain}
in Section~\ref{SSproof}.

In Part~\ref{PARTpolarizations} we prove a very general theorem that sheds
additional light on the proof of Theorem~\ref{TTmain}. In 
Section~\ref{SSattainablestatement} we associate to every isogeny class $\CC$ of
     \marginpar{\color{Red}\small\raggedright
                No --- the proof of this statement is incorrect.}
abelian varieties over a field $k$ a two-torsion group $\CBtwo(\CC)$ and a 
finite set $S_{\CC}\subseteq\CBtwo(\CC)$, and 
     {\color{Red} we prove in Section~\ref{SSattainableproof} 
      that the set $S_{\CC}$ determines the set of kernels of
      polarizations of varieties in $\CC$ up to Jordan-H\"older
      isomorphism.}
Then in Section~\ref{SSrevisit} we revisit the proof of Theorem~\ref{TTmain} and
show how it can be interpreted in terms of the group $\CBtwo(\CC)$ and the
set~$S_{\CC}$.

Theorem~\ref{TTmain} was inspired by a construction of Silverberg and Zarhin 
(see~\cite{brefSilverbergZarhinReps} and~\cite{brefSilverbergZarhinPolReps}); 
they too construct twists of a power of an elliptic curve such that every 
polarization of every abelian variety isogenous to the twist has degree
divisible by a given prime. Their original construction is limited to base 
fields of positive characteristic that have nonabelian Galois extensions of a
certain type, but more recently they have produced a new construction that works
over an arbitrary number field (see~\cite{brefSilverbergZarhinPols}).

\medskip\noindent
\textbf{Acknowledgments.}
The author thanks Daniel Goldstein, Bob Guralnick,\break 
Alice Silverberg, and Yuri Zarhin for helpful conversations and correspondence.

\medskip\noindent
\textbf{Conventions and notation.}
We consider varieties to be schemes over some specified base field; it follows
that if $U$ and $V$ are varieties over a field~$k$, then what we call a 
\emph{morphism} from $U$ to $V$ others might call a \emph{$k$-morphism} from
$U$ to~$V$. If $U$ is a variety over a field $k$ and $\ell$ is an extension 
field of~$k$, then we let $U_{\ell}$ denote the $\ell$-variety 
$V\times_{\Spec k} \Spec \ell$. If $\alpha\colon U\ra V$ is a morphism of
varieties over~$k$, we let $\alpha_{\ell}$ denote the induced morphism from 
$U_{\ell}$ to~$V_{\ell}$. If $X$ is an abelian variety, we let $\Xh$ denote its
dual variety, and if $\alpha\colon X\ra Y$ is a morphism of abelian varieties,
we let $\hat{\alpha}$ denote the dual morphism $\Yh\ra\Xh$. If~$G$ is a group
scheme over $k$ and $n$ is an integer, we denote by $G[n]$ the $n$-torsion 
subscheme of~$G$.

\section{Isogeny classes containing no principally polarized varieties}
\label{PARTnoPPAVs}

\subsection{Polarizations and endomorphisms of Galois twists of abelian varieties}
\label{SSgeneral}
In this section we prove some simple general results about Galois twists of 
abelian varieties that we will need in our proof of Theorem~\ref{TTmain}.

Suppose $k$ is a field, $\ell$ is a Galois extension of $k$ with Galois 
group~$G$, and $Y$ is an abelian variety over~$k$. Suppose that $X$ is an 
$\ell/k$-twist of $Y$ and that $f\colon Y_{\ell}\ra X_{\ell}$ is an isomorphism.
Then $X$ corresponds (as in Section~III.1.3 of~\cite{brefSerre}) to the element
of $H^{1}(G, \Aut Y_{\ell})$ represented by the cocycle 
$\sigma\mapsto a_{\sigma}:=f^{-1} f^{\sigma}$. Our first proposition tells us
when an endomorphism $\alpha$ of $Y$ gives rise to an endomorphism of~$X$.

\begin{proposition}
\label{PPendomorphism}
The endomorphism $f\alpha_{\ell} f^{-1}$ of $X_{\ell}$ descends to an 
endomorphism of $X$ if and only if we have 
$a_{\sigma}\alpha_{\ell} = \alpha_{\ell} a_{\sigma}$ for all $\sigma$ in~$G$.
\end{proposition}

\begin{proof}
Let $\beta$ be the endomorphism $f\alpha_{\ell} f^{-1}$ of~$X_{\ell}$. Then 
$\beta$ will descend to $X$ if and only if for all $\sigma$ in $G$ we have 
$\beta^{\sigma}=\beta$, which is 
\[ (f \alpha_{\ell} f^{-1})^{\sigma} = f \alpha_{\ell} f^{-1}. \]
By multiplying this equality by $f^{\sigma}$ on the right and by $f^{-1}$ on the
left, and by using the fact that $\alpha_{\ell} = \alpha_{\ell}^{\sigma}$, we
see that $\beta$ descends to $X$ if and only if for all $\sigma$ we have
\[ f^{-1} f^{\sigma} \alpha_{\ell} = \alpha_{\ell}  f^{-1} f^{\sigma}, \]
if and only if $a_{\sigma}\alpha_{\ell} = \alpha_{\ell} a_{\sigma}$ for all
$\sigma$ in~$G$.
\end{proof}

Now suppose $\lambda$ is a polarization of~$Y$, and let $x\mapsto x^{\dagger}$ 
denote the Rosati involution on $\End Y$ corresponding to~$\lambda$, so that 
$x^{\dagger} = \lambda^{-1}{}\hat{x}{}\lambda$. Our second proposition tells us 
when the polarization $\lambda$ gives rise to a polarization of~$X$.

\begin{proposition}
\label{PPpolarization}
The polarization $\hat{f^{-1}} \lambda_{\ell} f^{-1}$ of $X_{\ell}$ descends to
a polarization of $X$ if and only if we have 
$a_{\sigma}^{\dagger} a_{\sigma} = 1$ for all $\sigma$ in~$G$.
\end{proposition}

\begin{proof}
Let $\mu$ be the polarization $\hat{f^{-1}} \lambda_{\ell} f^{-1}$ 
of~$X_{\ell}$. Then $\mu$ will descend to $X$ if and only if for all $\sigma$ in
$G$ we have $\mu^{\sigma}=\mu$, which is
\[(\hat{f^{-1}} \lambda_{\ell}  f^{-1})^{\sigma} 
     = \hat{f^{-1}} \lambda_{\ell} f^{-1}.\]
By multiplying this equality by $f^{\sigma}$ on the right and by 
$\lambda_{\ell}^{-1} \hat{f^{\sigma}}$ on the left, and by using the fact that
$\lambda_{\ell} = \lambda_{\ell}^{\sigma}$, we see that $\mu$ descends to $X$ if
and only if for all $\sigma$ we have
\[\lambda_{\ell}^{-1}  \hat{(f^{-1} f^{\sigma})}\lambda_{\ell}
     (f^{-1} f^{\sigma}) = 1,\]
if and only if $a_{\sigma}^{\dagger} a_{\sigma} = 1$ for all $\sigma$ in~$G$.
\end{proof}

\vspace{-0.8ex} 

\subsection{Reduced restrictions of scalars of $p$-isolated elliptic curves}
\label{SSrestriction}
Let $k$ be a $p$-admissible field, let $E$ be a $p$-isolated elliptic curve 
over~$k$, and let $\ell$ be a Galois extension of $k$ of degree~$p$. Let $X$ be
the reduced restriction of scalars of $E$ from $\ell$ to~$k$. In this section we
calculate the endomorphism ring of~$X$, a restriction on the degrees of the
polarizations of~$X$, and the Galois module structure of the $p$-torsion of~$X$.
These results are the building blocks of our proof of Theorem~\ref{TTmain}.

\begin{lemma}
\label{LLnoCM}
The elliptic curve $E_{\ell}$ has no complex multiplication.
\end{lemma}

\begin{proof}
The endomorphism algebra $A = (\End E_\ell)\otimes\BQ$ is either~$\BQ$, an
imaginary quadratic field, or a quaternion algebra over~$\BQ$. The action of 
$\Gal(\ell/k)$ on $\End E_\ell$ gives us a homomorphism 
$\Gal(\ell/k)\ra \Aut(A/\BQ)$, and since $\End E = \BZ$, this homomorphism must
be nontrivial if $A\neq\BQ$. Suppose $A$ were a quaternion algebra. Then the
image of a generator of $\Gal(\ell/k)$ in $\Aut(A/\BQ)$ would be a nontrivial
automorphism. Since every automorphism of a quaternion algebra is inner, this
automorphism would have to be given by conjugation by a non-central element $s$
of~$A$. But then the $2$-dimensional sub-algebra $\BQ(s)$ would be fixed by the
action of $\Gal(\ell/k)$ on~$A$, contradicting the fact that $\End E = \BZ$. 
Therefore $A$ is not a quaternion algebra. Suppose $A$ were a quadratic field. 
Then $\Aut(A/\BQ)$ would be a cyclic group of order~$2$, contradicting the 
existence of a nontrivial homomorphism $\Gal(\ell/k)\ra \Aut(A/\BQ)$. Thus $A$
must be~$\BQ$, and $\End E_{\ell}$ must be~$\BZ$.
\end{proof}

Let $\sigma$ be a generator of $\Gal(\ell/k)$ and let $R$ be the restriction of
scalars of $E$ from $\ell$ to~$k$. Then $R$ is the $\ell/k$-twist of $E^{p}$
given by the element of $H^{1}(\Gal(\ell/k), \Aut E_{\ell}^{p})$ represented by
the cocycle that sends $\sigma$ to the automorphism $\xi$ of $E_{\ell}^{p}$ that
cyclically shifts the factors. The kernel $S$ of the trace map $E^{p}\ra E$ is
stable under~$\xi$, and the reduced restriction of scalars $X$ of $E$ is the 
$\ell/k$-twist of $S$ given by sending $\sigma$ to the restriction of $\xi$ 
to~$S$. The projection map from $E^{p}$ onto its first $p-1$ factors gives an
isomorphism from $S$ to $E^{p-1}$; under this isomorphism, the restriction of
$\xi$ to $S$ is given by the $(p-1)\times(p-1)$ matrix \nobreak
\[
\zeta = \left[\begin{matrix}
-1                & -1                & \cdots\ & \llap{$-$}1     & -1     \\
\phantom{-}1      & \phantom{-}0      & \cdots\ &          0      & \phantom{-}0      \\
\phantom{-}0      & \phantom{-}1      & \cdots\ &          0      & \phantom{-}0      \\
\phantom{-}\vdots & \phantom{-}\vdots & \ddots\ &          \vdots & \phantom{-}\vdots \\
\phantom{-}0      & \phantom{-}0      & \cdots\ &          1      & \phantom{-}0      \\
\end{matrix} \ \right],
\]
where we identify the ring $M_{p-1}(\BZ)$ of $(p-1)\times(p-1)$ integer matrices
with the endomorphism ring of $E_{\ell}^{p-1}$. In other words, $X$ is the
$\ell/k$-twist of $E^{p-1}$ given by the element of
$H^{1}(\Gal(\ell/k), \Aut E_{\ell}^{p-1})$ represented by the cocycle that sends
$\sigma$ to~$\zeta$. Note that the minimal polynomial of the endomorphism 
$\zeta$ is the $p$th cyclotomic polynomial.

\begin{lemma}
\label{LLendomorphism}
The abelian variety $X$ is simple over~$k$, and its endomorphism ring is
isomorphic to the ring of integers of the $p$th cyclotomic field.
\end{lemma}

\begin{proof}
Since $X_{\ell} \cong E_{\ell}^{p-1}$ and $\End E_{\ell} = \BZ = \End E$, every
endomorphism of $X$ comes from an endomorphism of~$E^{p-1}$. According to
Proposition~\ref{PPendomorphism}, the only endomorphisms of $E^{p-1}$ that give
rise to elements of $\End X$ are the endomorphisms that commute with the element
$\zeta$ of $\End E_{\ell}^{p-1}$ defined by the matrix above. Since $\BQ(\zeta)$
is a field of degree $p-1$ over~$\BQ$, it is a maximal commutative subring of
the matrix ring $M_{p-1}(\BQ)$, so the only elements of $\End E^{p-1}$ that
commute with $\zeta$ are those elements that lie in~$\BQ(\zeta)$. The 
intersection of $\BQ(\zeta)$ with $\End E^{p-1}$ is $\BZ[\zeta]$, which is the
ring of integers of the cyclotomic field~$\BQ(\zeta)$. Thus $\End X$ is 
isomorphic to the ring of integers of the $p$th cyclotomic field. And finally,
the fact that $(\End X)\otimes\BQ$ is a field shows that $X$ is simple.
\end{proof}

\begin{lemma}
\label{LLdegree}
If $\alpha\in\BQ(\zeta)$ is an endomorphism of~$X$, then the degree of $\alpha$
is the square of the norm of $\alpha$ from $\BQ(\zeta)$ to~$\BQ$.
\end{lemma}

\begin{proof}
It is easy to see that under the identification of $\End E^{p-1}$ with
$M_{p-1}(\BZ)$, the degree function is the square of the determinant function.
The lemma then follows from the fact the determinant function from 
$M_{p-1}(\BQ)$ to~$\BQ$, restricted to the maximal subfield $\BQ(\zeta)$ of 
$M_{p-1}(\BQ)$, is the field norm from $\BQ(\zeta)$ to~$\BQ$.
\end{proof}

\begin{lemma}
\label{LLonepolarization}
The abelian variety $X$ has a polarization $\lambda$ of degree~$p^{2}$.
\end{lemma}

\begin{proof}
Let $\mu$ be the product principal polarization of~$E^{p-1}$, let $b$ be the
endomorphism of $E^{p-1}$ defined by the matrix
\[
\left[\begin{matrix}
2      & 1      & \cdots & 1      \\
1      & 2      & \cdots & 1      \\
\vdots & \vdots & \ddots & \vdots \\
1      & 1      & \cdots & 2      \\
\end{matrix}\right]
\]
with $2$'s on the diagonal and $1$'s elsewhere, and let $\bar{\lambda} = \mu b$.
The Rosati involution associated to $\mu$ is the matrix transpose 
$x\mapsto~x^{t}$, and the Rosati involution associated to $\bar{\lambda}$ is the
matrix transpose conjugated by $b$ --- that is, $x\mapsto b^{-1} x^{t} b$. One
checks that $b^{-1} \zeta^{t} b \zeta$ is the identity matrix, so by 
Proposition~\ref{PPpolarization} the polarization ${\bar{\lambda}}_{\ell}$ of
$E_{\ell}^{p-1}$ descends to give a polarization $\lambda$ of~$X$.

\goodbreak

Since $\mu$ is a principal polarization and 
$\lambda_{\ell} = \mu_{\ell} b_{\ell}$, the degree of $\lambda$ is the degree
of~$b$, which is the square of the determinant of~$b$. An easy calculation shows
that $\det b = p$, which proves the lemma.
\end{proof}

\begin{lemma}
\label{LLallpolarizations}
Suppose $\mu$ is a polarization of~$X$. Then there is an integer $n$\break
such that $\deg\mu = p^{2}n^{4}$.
\end{lemma}

\begin{proof}
Let $K$ be the field $(\End X)\otimes\BQ$ and let $\zeta$ be the endomorphism of
$X$ defined above, so that $\zeta$ is a  primitive $p$th root of unity and 
$K = \BQ(\zeta)$. Let $x\mapsto x^{\dagger}$ be the Rosati involution associated
to the polarization $\lambda$ of Lemma~\ref{LLonepolarization}. Then for every
$x\in K$, the element $x^{\dagger}$ is simply the complex conjugate of~$x$.

Suppose $\mu$ is a polarization of~$X$. Then the element $\lambda^{-1}\mu$ of
the field $K$ is fixed by the Rosati involution (see \S21, Application~3,
pp.~208--210 of~\cite{brefMumford}), and is therefore an element of the maximal
real subfield $K^{+}$ of~$K$. By Lemma~\ref{LLdegree}, the degree of an element
of $K$ is equal to the square of its norm to~$\BQ$, so we find that
\[
\deg\mu = \deg\lambda\cdot\deg(\lambda^{-1}\mu) 
= p^{2} \left(N_{K/\BQ}(\lambda^{-1}\mu)\right)^{2}
= p^{2} \left(N_{K^{+}/\BQ}(\lambda^{-1}\mu)\right)^{4}.
\]
Let $n = N_{K^{+}/\BQ}(\lambda^{-1}\mu)$. Since $\deg\mu$ is an integer and $n$
is rational, we see that $n$ is an integer, and the lemma is proved.
\end{proof}

We end this section by describing the Galois module structure of~$X[p]$. First 
note that $p$ is not equal to the characteristic of the base field $k$ because
$E$ has no $p$-isogenies. Thus, the group scheme structures of $E[p]$ and $X[p]$
are completely captured by the Galois module structures of the sets of points of
these schemes over a separable closure $\kb$ of~$k$. Furthermore, $E[p](\kb)$ is
a simple $\Gal(\kb/k)$-module because $E$ has no $p$-isogenies to another
elliptic curve over~$k$.

\begin{lemma}
\label{LLgaloismodule}
The sequence of modules
\[ 
0 \subset (\zeta-1)^{p-2} X[p](\kb) 
  \subset (\zeta-1)^{p-3} X[p](\kb) 
  \subset \cdots 
  \subset X[p](\kb)
\]  
is a composition series for the $\Gal(\kb/k)$-module~$X[p](\kb)$, and each
composition factor is isomorphic to~$E[p](\kb)$.
\end{lemma}

\begin{proof}
Multiplication by $(\zeta-1)^{p-i-2}$ gives an isomorphism from the quotient
$(\zeta-1)^{i}X[p](\kb) / (\zeta-1)^{i+1}X[p](\kb)$ to 
$(\zeta-1)^{p-2} X[p](\kb)$, which is the kernel of $\zeta-1$ acting 
on~$X[p](\kb)$. This kernel is the image of $E[p](\kb)$ under the diagonal 
embedding of $E_{\ell}$ into~$E_{\ell}^{p-1}$.
\end{proof}

\begin{lemma}
\label{LLisogenymodule}
If $\varphi\colon X\ra Y$ is an isogeny, then the only simple 
$\Gal(\kb/k)$-module that occurs as a composition factor of the 
$p$-power-torsion part of the kernel of $\varphi$ is~$E[p](\kb)$.
\end{lemma}

\begin{proof}
Immediate from Lemma~\ref{LLgaloismodule}.
\end{proof}

\subsection{Proof of Theorem~\ref{TTmain}}
\label{SSproof}
For every finite group scheme $G$ over~$k$, let us define the \emph{$E[p]$-rank}
of $G$ to be the multiplicity of the simple Galois module $E[p](\kb)$ as a
composition factor of the $p$-power-torsion of~$G$. We will denote the 
$E[p]$-rank of $G$ by $\rank_{E[p]}(G)$. The $E[p]$-rank is an additive function
on exact sequences of finite group schemes.

Suppose $\varphi\colon X\ra Y$ is an isogeny and $\nu$ is a polarization of~$Y$.
Then the map $\mu \colon X\ra \Xh$ given by $\hat{\varphi} \nu \varphi$ is a
polarization of~$X$, so we have
\[
\rank_{E[p]}(\ker\mu) 
   =   \rank_{E[p]}(\ker\hat{\varphi})
     + \rank_{E[p]}(\ker\nu) 
     + \rank_{E[p]}(\ker\varphi).
\]
Now, $\ker\hat{\varphi}$ is the Cartier dual of $\ker\varphi$, and $E[p]$ is its
own Cartier dual, so 
$\rank_{E[p]}(\ker\hat{\varphi}) = \rank_{E[p]}(\ker\varphi)$. Thus the parity 
of $\rank_{E[p]}(\ker\nu)$ is equal to that of $\rank_{E[p]}(\ker\mu)$, which is
odd by Lemmas~\ref{LLallpolarizations} and~\ref{LLisogenymodule}. Therefore
$E[p]$ appears as a composition factor of the $p$-power-torsion of $\ker\nu$, so 
$p^{2}$ divides the order of $\ker\nu$, so $p^{2}$ divides the degree of~$\nu$.
\qed

\section{Polarizations up to Jordan-H\"older isomorphism}
\label{PARTpolarizations}

\noindent
In this part of the paper we associate to every isogeny class $\CC$ of abelian 
varieties over a field $k$ a two-torsion group $\CBtwo(\CC)$ and a finite set 
$S_{\CC}\subseteq\CBtwo(\CC)$, and we show that the set $S_{\CC}$ determines the
set of kernels of polarizations of varieties in $\CC$ up to Jordan-H\"older 
isomorphism. Then we revisit the proof of Theorem~\ref{TTmain} and show how it
can be interpreted in terms of the group $\CBtwo(\CC)$ and the set~$S_{\CC}$.

\subsection{Statement of results}
\label{SSattainablestatement}
For every isogeny class $\CC$ of abelian varieties over a field~$k$, we let
$\Ker_\CC$ be the category whose objects are finite commutative group schemes 
over $k$ that can be embedded (as closed sub-group-schemes) in some variety in
the isogeny class $\CC$ and whose morphisms are morphisms of group schemes. We 
see that the objects in $\Ker_\CC$ are those group schemes that can be written 
$\ker\varphi$ for some isogeny $\varphi\colon X\ra Y$ of elements of~$\CC$.

The \emph{Grothendieck group} $G(\Ker_{\CC})$ of $\Ker_{\CC}$ is defined to be
the quotient of the free abelian group on the isomorphism classes of objects in
$\Ker_{\CC}$ by the subgroup generated by the expressions $X-X'-X''$ for all
exact sequences $0\ra X'\ra X\ra X''\ra 0$ in $\Ker_{\CC}$. If $X$ is an object 
of $\Ker_\CC$ we let $[X]$ denote its class in $G(\Ker_\CC)$, and we say that
two objects $X$ and $Y$ are \emph{Jordan-H\"older isomorphic} to one another if
$[X] = [Y]$. The group $G(\Ker_\CC)$ is a  free abelian group on the simple 
objects of $\Ker_{\CC}$. An element of $G(\Ker_\CC)$ is said to be 
\emph{effective} if it is a sum of positive multiples of simple objects. Let us
call an element $P$ of $G(\Ker_{\CC})$ \emph{attainable} if there is a
polarization $\lambda$ of a variety in $\CC$ such that $P = [\ker\lambda]$. Our
goal in this section will be to identify the attainable elements of 
$G(\Ker_{\CC})$.

To identify the attainable elements we must first define several groups for
every isogeny class of abelian varieties. The first two groups will be subgroups
of $G(\Ker_{\CC})$, and the others will be defined solely in terms of the
endomorphism rings of the varieties in the isogeny class.

Let $Z(\CC)$ denote the subgroup of $G(\Ker_{\CC})$ generated by the elements of
the form~$[G]$, where $G \in \Ker_{\CC}$ is a group scheme whose rank is a 
square and for which there exists a non-degenerate alternating pairing 
$G\times G\ra\BGm$.\footnote{
     The existence of a non-degenerate alternating pairing \emph{implies}
     that the rank of $G$ is a square, except in characteristic~$2$; the
     unique simple local-local group scheme in characteristic~$2$, which
     has rank~$2$, has a non-degenerate alternating pairing.}
Cartier duality on $\Ker_{\CC}$ defines an involution $P\mapsto\overline{P}$ of
$G(\Ker_{\CC})$, and $Z(\CC)$ is stable under this involution. Let $B(\CC)$ be 
the subgroup $\{ P + \overline{P} : P \in G(\Ker_{\CC})\}$ of $G(\Ker_{\CC})$; 
it is not hard to see that $B(\CC) \subseteq Z(\CC)$.

Now let us define the groups that depend on the endomorphism rings of the
varieties in~$\CC$. If $X$ and $Y$ are two varieties in $\CC$ then 
$(\End X)\otimes\BQ \cong (\End Y)\otimes\BQ$, so we may define $\End^{0} \CC$
to be $(\End X)\otimes\BQ$ for any $X$ in~$\CC$. We may write 
$\CC \sim  \CC_{1}^{n_{1}}\times \cdots \times \CC_{r}^{n_{r}}$ for some 
distinct isogeny classes $\CC_{i}$ of simple varieties and some integers 
$n_{i}$; by this we mean that every $X$ in $\CC$ is isogenous to a product 
$X_{1}^{n_{1}}\times \cdots \times X_{r}^{n_{r}}$ where $X_{i}\in \CC_{i}$. Let 
$A = \End^{0}\CC$. Then 
\[
A \cong M_{n_{1}}(D_{1}) \times \cdots\times M_{n_{r}}(D_{r}),
\]
where $D_{i} = \End^{0}\CC_{i}$ and where $M_{n_{i}}(D_{i})$ denotes an 
$n_{i}$-by-$n_{i}$ matrix algebra over~$D_{i}$. Each $D_{i}$ is a division 
algebra over its center~$K_{i}$, which is a number field, and the center $K$ of
$A$ is the product of the~$K_{i}$. Only certain kinds of division algebras can 
occur as the endomorphism algebras of simple isogeny classes, and they are 
classified into four types (see Theorem~2 (p.~201) of~\cite{brefMumford}). We
will define three subgroups $R_{0}(A)$, $R_{1}(A)$, and $R_{2}(A)$ of $K^{*}$ by
defining these groups first for the four possible types of division algebras $D$ 
and by then setting
\[
R_{i}(A) = R_{i}(D_{1}) \times \cdots\times R_{i}(D_{r}).
\]
The group $R_{1}(A)$ will be a subgroup of $R_{0}(A)$, and $R_{2}(A)$ will be a
subgroup of finite index in $R_{1}(A)$.

\textbf{Type I:} For this type, $D = K$ is a totally real number field. We let
$R_{0}(D) = K^{*}$, we let $R_{1}(D)$ be the multiplicative group of totally 
positive elements of~$K$, and we let $R_{2}(D) = R_{1}(D)$.

\textbf{Type II:} For this type, $K$ is a totally real number field and $D$ is a
quaternion algebra over $K$ that is split at every infinite prime of~$K$. We let
$R_{0}(D) = K^{*}$, we let $R_{1}(D)$ be the multiplicative group of totally
positive elements of~$K$, and we let $R_{2}(D)$ be the subgroup of $R_{1}(D)$
consisting of those elements that are squares in $K_{\pid}$ for all primes 
$\pid$ of $K$ for which the local Brauer invariant $\inv_{\pid} D$ of $D$ at
$\pid$ is nonzero.

\textbf{Type III:} For this type, $K$ is a totally real number field and $D$ is a 
\hbox to \hsize{quaternion algebra over $K$ that is ramified at everyinfinite prime of~$K$. We}

\noindent
let $R_{0}(D)$ be the multiplicative group of totally positive elements of~$K$,
we let $R_{1}(D)$ be the group of squares of elements of $R_{0}(D)$, and we let
$R_{2}(D) = R_{1}(D)$.

\textbf{Type IV:} For this type, $K$ is a CM-field with maximal real subfield
$K^{+}$ and $D$ is a division algebra over $K$ such that if $\sigma$ is the 
nontrivial automorphism of $K$ over $K^{+}$ then $\inv_{\pid} D = 0$ for every
prime of $K$ fixed by $\sigma$ and $\inv_{\pid} D + \inv_{\pid^{\sigma}} D = 0$
for every prime $\pid$ of~$K$. We let $R_{0}(D) = K^{*}$, we let $R_{1}(D)$ be
the multiplicative group of totally positive elements of~$K^{+}$, and we let
$R_{2}(D) = R_{1}(D)$.

An involution $x\mapsto x^{\dagger}$ of $A$ is \emph{positive} if 
$\Trd(x x^{\dagger})$ is a totally positive element of $K$ for every element $x$
of~$A^{*}$, where $\Trd$ is the reduced trace map from $A$ to~$K$. If 
$\alpha\in A$ is fixed by a positive involution then $\BQ(\alpha)$ is a product
of totally real number fields. If $x\mapsto x^{\dagger}$ is a positive 
involution of~$A$, we let $A_{\dagger}$ denote the set of elements of $A^{*}$
that are fixed by $\dagger$ and that are totally positive.

The following lemma, whose proof we will give in 
Section~\ref{SSattainableproof}, motivates the definitions of the groups 
$R_{i}(A).$

\begin{lemma}
\label{LLRexplain}
Let $A = \End^{0}\CC$ for an isogeny class $\CC$ of abelian varieties over~$k$.
Then $R_{0}(A) = \Nrd(A^{*})$, where $\Nrd$ is the reduced norm map from $A$ to
its center~$K$. Suppose $x\mapsto x^{\dagger}$ is a positive involution of~$A$. 
Then $R_{2}(A)\subseteq \Nrd(A_{\dagger}) \subseteq R_{1}(A).$
\end{lemma}

\noindent
\textbf{Remark.}
Note that if $A$ is built up out of simple $D$ that are of type I, III, and~IV
then $R_{2}(A) = R_{1}(A)$ and $\Nrd(A_{\dagger})$ is a subgroup of~$K^{*}$, but
if a $D$ of type~II occurs in $A$ then Lemma~\ref{LLRexplain} only allows us to
say that $\Nrd(A_{\dagger})$ lies between two subgroups of~$K^{*}$. We cannot
expect to do much better than this; one can find examples of $D$ of type II for
which $\Nrd(D_{\dagger})$ is not a group.

\medskip

Let $A = \End^{0}\CC$. We define a homomorphism $\Prin$ from the multiplicative 
group $A^{*}$ to $G(\Ker_{\CC})$ as follows: Suppose $\alpha\in A^{*}$ is given.
We pick a variety $X$ in $\CC$ and choose an endomorphism $\beta$ of $X$ and an
integer $n$ such that $\alpha = \beta /n$. Then we set 
$\Prin(\alpha) = [\ker\beta] - [\ker n]$. We leave it to the reader to show that
the value $[\ker\beta] - [\ker n]$ does not depend on the choice of~$X$,
$\beta$, and~$n$. 

Note that Lemma~\ref{LLRexplain} states in part that $R_{0}(A) = \Nrd(A^{*})$. 
We can use this identity, together with the homomorphism $\Prin$, to define a 
homomorphism $\Phi$ from $R_{0}(A)$ to $G(\Ker_{\CC})$. Suppose $a\in R_{0}(A)$
is given, and suppose $\alpha$ and $\alpha'$ are elements of $A^{*}$ with 
$\Nrd(\alpha) = \Nrd(\alpha') = a$. Then, by an easy consequence of Wang's 
Theorem (see~\cite{brefPlatonovRapinchuk}, Theorem~1.14, p.~38, and the comments
on p.~39), we see that $\alpha'\alpha^{-1}$ lies in the commutator subgroup
of~$A^{*}$, so that $\Prin(\alpha'\alpha^{-1}) = 0$. Thus we may define $\Phi$
by taking $\Phi(a) = \Prin(\alpha)$ for any choice of $\alpha\in A^{*}$ such
that $\Nrd(\alpha) = a$.

Finally, we let $\CBone(\CC)$ be the quotient of $Z(\CC)$ by the subgroup 
generated by 
\hbox to \hsize{$B(\CC)$ and $\Phi(R_{1}(A))$, and we let $\CBtwo(\CC)$ be the
quotient of $Z(\CC)$ by the subgroup}

\noindent
generated by $B(\CC)$ and $\Phi(R_{2}(A))$. Since $2 Z(\CC)$ is contained in 
$B(\CC)$ the groups $\CBone(\CC)$ and $\CBtwo(\CC)$ are $2$-torsion, and since
$R_{2}(A)$ is a subgroup of finite index in $R_{1}(A)$, the natural surjection
$\CBtwo(\CC)\ra\CBone(\CC)$ has a finite kernel.

\begin{theorem}
\label{TTattainable}
     \marginpar{\color{Red}\small\raggedright
                This result is incorrect. With the definition of 
                $S_{\CC}$ given below, the ``only if'' implication
                holds but the ``if'' implication does not. The error
                is indicated below.}
There is a finite subset $S_{\CC}$ of the group $\CBtwo(\CC)$ such that an 
element of $G(\Ker_{\CC})$ is attainable 
     {\color{Red} if and only if}
it is an effective element of $Z(\CC)$ whose image in $\CBtwo(\CC)$ lies
in~$S_{\CC}$. Furthermore, the image of $S_{\CC}$ in $\CBone(\CC)$ consists of a
single element~$I_{\CC}$.
\end{theorem}

This theorem must surely seem quite removed from the down-to-earth question of 
whether there exists a principally-polarized variety in a given isogeny class.
However, we will show in Section~\ref{SSrevisit} that our proof of 
Theorem~\ref{TTmain} can be viewed as an argument showing that there is a 
homomorphism $\CBone(\CC)\ra\BZ/2\BZ$ such that the image of $I_{\CC}$ is 
nonzero. Also, in certain cases the abstract objects in 
Theorem~\ref{TTattainable} can be computed. For example, if $\CC$ is an isogeny
class over a finite field then the group $\CBone(\CC)$ (which in this case is
equal to $\CBtwo(\CC)$) can be computed, and one can even compute a subgroup of
order at most $2$ that contains $I_{\CC}$ (see \cite{brefHoweAV}). If in 
addition the varieties in $\CC$ are ordinary, the element $I_{\CC}$ itself can 
be calculated (see \cite{brefHoweOAV}).

\subsection{Proofs of Lemma~\ref{LLRexplain} and Theorem~\ref{TTattainable}}
\label{SSattainableproof}
In this section we prove the results stated in the preceding section.

\begin{proof}[Proof of Lemma~{\textup{\ref{LLRexplain}}}]
Clearly it will suffice to prove the lemma in the case where $A=M_{n}(D)$ for
some integer $n$ and division algebra $D$ of one of the four types listed in the
preceding section. For each of these types of algebras, the 
Hasse-Schilling-Maass theorem (\cite{brefReiner}, Theorem~33.15, p.~289) shows
that $R_{0}(A) = \Nrd(A^{*})$, so we need only prove the second statement of the
lemma.

If $D/K$ is of Type III then the statement we are to prove is Theorem~4.7 
of~\cite{brefHoweAV}, while if $D/K$ is of Type IV then the statement we are to
prove is Theorem~4.1 of~\cite{brefHoweAV}.

Suppose $D=K$ is of Type I. We must show that $\Nrd(A_{\dagger})$ is the set of
totally positive elements of~$K$. It is clear that $\Nrd(\alpha)$ is totally
positive if $\alpha\in A_{\dagger}$, so all we must show is that every totally
positive element of $K$ is the reduced norm of some element of~$A_{\dagger}$.

Let $x\mapsto x^*$ be the transpose on~$M_d(K)$. Then by Theorem~8.7.4
(pp.~301--302) and Theorem~7.6.3 (p.~259) of~\cite{brefScharlau}, there is an
isomorphism $i\colon A\ra M_d(K)$ and a diagonal matrix $\alpha\in M_d(K)$ such
that the isomorphism $i$ takes the involution $x\mapsto x^{\dagger}$ to the
involution $\eta$ of $M_d(K)$ defined by $\eta(x)=\alpha x^* \alpha^{-1}$. Now
suppose we are given a totally positive $b$ in~$K$. Let $\beta\in M_{d}(K)$ be
the diagonal matrix with $b$ in the upper left corner and $1$'s elsewhere. Then
$\beta$ is totally positive and fixed by~$\eta$, and its determinant is~$b$.
Therefore $i^{-1}(\beta)$ is an element of $K_{\dagger}$ with reduced norm equal
to~$b$. 

\hbox to \hsize{\hskip\parindent Suppose $D/K$ is of Type II. It is clear that the reduced norm of an element}

\noindent
of $A_{\dagger}$ is totally positive, so we have 
$\Nrd(A_{\dagger}) \subseteq R_{1}(A)$, and we must prove that 
$R_{2}(A)\subseteq \Nrd(A_{\dagger})$.

Let $x\mapsto x^*$ be the conjugate transpose involution on $M_d(D)$, where 
``conjugation'' on $D$ is the standard involution $x\mapsto \Trd_{D/K} x - x$. 
Then again by Theorem~8.7.4 (pp. 301--302) and Theorem~7.6.3 (p.~259) 
of~\cite{brefScharlau}, there is an isomorphism $i\colon A\ra M_d(D)$ and a
diagonal matrix $\alpha\in M_d(D)$ with $\alpha^{*} = -\alpha$ such that the 
isomorphism $i$ takes the involution $x\mapsto x^{\dagger}$ to the involution 
$\eta$ of $M_d(D)$ defined by $\eta(x)=\alpha x^* \alpha^{-1}$; furthermore, as
is argued on pp.~194--195 of~\cite{brefMumford}, the entries of the diagonal
matrix $\alpha^{2}$ are totally negative elements of~$K$. Let 
$\alpha_{1},\ldots,\alpha_{n}$ denote the diagonal entries of $\alpha$ and let 
$c_{1} = \alpha_{1}^{2}$, so that $\alpha_{1}$ satisfies the polynomial 
$x^{2} - c_{1}$. Since the field $K(\alpha_{1})$ splits the quaternion 
algebra~$D$, the element $c_{1}$ of $K$ must be a nonsquare in $K_{\pid}$ for
every prime $\pid$ of $K$ for which $\inv_{\pid} D\neq 0$.

Suppose $b$ is an element of $R_{2}(A)$. We claim that there exists an element
$\beta$ of $D$ such that (1) $\beta + \beta^{*} = 0$ and (2) 
$\beta\beta^{*} = -bc_{1}$ and (3) $\beta\alpha_{1}^{-1}$ is totally positive. 
To see that such a $\beta$ exists, we first note that $bc_{1}$ is a nonsquare in
$K_{\pid}$ for every prime $\pid$ of $K$ for which $\inv_{\pid} D\neq 0$, so the
field $K(\sqrt{bc_{1}})$ splits~$D$. This shows that there is an element 
$\beta_{0}$ of $D$ such that $\beta_{0}^{2}  - bc_{1} = 0$, so that $\beta_{0}$ 
satisfies (1) and (2). If we choose a $K$-basis for the trace-$0$ elements
of~$D$, then the set of $\beta$ satisfying (1) and (2) is a level set of a 
homogeneous ternary quadratic form $Q$ that is indefinite at every infinite
prime of~$K$, and we have just shown that there are $K$-points in this level
set. Condition (3) is simply a linear inequality at each of the infinite primes
of~$K$, and the inequality can be satisfied locally at each infinite prime by
points on the level set because the form $Q$ is indefinite, so by weak 
approximation we see that there do exist $\beta$'s satisfying (1), (2), and (3).

Choose such a $\beta\in D$, and let $\gamma$ be the diagonal matrix with 
$\beta \alpha_{1}^{-1}$ in the upper left corner and $1$'s elsewhere on the
diagonal. A computation shows that $\gamma$ is fixed by the involution~$\eta$,
that the reduced norm of $\gamma$ is 
$\Nrd_{D/K}(\beta)/\Nrd_{D/K}(\alpha_{1}) = b$, and that $\gamma$ is totally 
positive. Thus $i^{-1}(\gamma)$ is an element of $A_{\dagger}$ with reduced norm
equal to~$b$.
\end{proof}

\begin{proof}[Proof of Theorem~{\textup{\ref{TTattainable}}}]
Suppose $\lambda\colon X\ra\Xh$ and $\mu\colon Y\ra\Yh$ are polarizations of
varieties in~$\CC$. First note that the ranks of $\ker\lambda$ and $\ker\nu$ are
squares, and that there are non-degenerate alternating pairings from these 
groups to the multiplicative group (see~\S23 of~\cite{brefMumford}), so that
$[\ker\lambda]$ and $[\ker\mu]$ lie in~$Z(\CC)$. Let $\varphi\colon X\ra Y$ be 
an isogeny, and let $\nu$ be the polarization $\hat{\varphi}\mu\varphi$ of~$X$,
where $\hat{\varphi}\colon\Yh\ra\Xh$ is the dual isogeny of~$\varphi$. Let $n$
be any positive integer such that $\ker\lambda\subseteq\ker(n\nu)$ as group 
schemes. Then there is an isogeny $\alpha\colon\Xh\ra\Xh$ such that 
$n\nu=\alpha\lambda$. A polarization is equal to its own dual isogeny, so we can
equate the right-hand side of the last equality with its dual to get 
$n\nu=\lambda\hat{\alpha}$. Using Application~III (pp.~208--210) of~\S21
of~\cite{brefMumford} (see especially the final paragraph) and the fact that 
$n\nu$ and $\lambda$ are polarizations, we find that $\hat{\alpha}\in \End X$ is
fixed by the Rosati involution associated to $\lambda$ and is totally positive.

The equality $n\hat{\varphi}\mu\varphi=\lambda\hat{\alpha}$ translates into the 
equality 
\[
[\ker n]+[\ker\hat{\varphi}]+[\ker\mu]+[\ker\varphi] 
  = [\ker\lambda]+[\ker\hat{\alpha}]
\]
in $G(\Ker_\CC)$. Now, $[\ker\hat{\varphi}]+[\ker\varphi]$ is an element of 
$B(\CC)$, and $\hat{\alpha}$ and $n$ are totally positive elements of 
$\End^{0}\CC$ that are fixed by the Rosati involution, so $[\ker n]$ and 
$[\ker\hat{\alpha}]$ lie in $\Phi(R_{1}(\End^{0}\CC))$. It follows that the 
images of $[\ker\mu]$ and $[\ker\lambda]$ in $\CBone(\CC)$ are equal. Thus, we 
may define $I_{\CC}$ to be the image in $\CBone(\CC)$ of the kernel of any 
polarization of any variety in~$\CC$.

Let $S_{\CC}$ be the image in $\CBtwo(\CC)$ of the subset 
     \marginpar{\color{Red}\small\raggedright 
                Note that this choice of $S_{\CC}$ is the only one
                that will make the ``only if'' implication of 
                Theorem~\ref{TTattainable} hold.}
\[
\{[\ker\lambda] : 
  \text{$\lambda$ is a polarization of some $X$ in $\CC$}\}
\]  
of $G(\Ker_{\CC})$. We see that $S_{\CC}$ is a subset of the preimage of 
$I_{\CC}$ under the natural reduction map $\CBtwo(\CC)\ra\CBone(\CC)$. Since 
this reduction map has a finite kernel, the set $S_{\CC}$ is finite.

To complete the proof of Theorem~\ref{TTattainable}, we must show that if $P$ is
an effective element of $G(\Ker_{\CC})$ whose image in $\CBtwo(\CC)$ lies 
in~$S_{\CC}$, then $P$ is attainable. So suppose $P$ is an effective element of 
$G(\Ker_{\CC})$ whose image in $\CBtwo(\CC)$ is equal to the image of 
$[\ker\lambda]$ for a polarization $\lambda\colon X\ra\Xh$ of a variety 
in~$\CC$. Then there is an element $Q$ of $G(\Ker_{\CC})$ and an element $a$ of
$R_{2}(\End^{0}\CC)$ such that $P + Q + \Qb = [\ker\lambda] + \Phi(a)$ in 
$G(\Ker_{\CC})$. Lemma~\ref{LLRexplain} shows that there is an 
$\alpha\in\End^{0}\CC$ that is totally positive and fixed by the Rosati 
involution associated to $\lambda$ such that $\Nrd\alpha = a$. Choose an integer
$n$ so that $n^{2}\alpha$ is an actual endomorphism of $X$ and such that 
$Q  + [\ker n]$ is effective. Replacing $\alpha$ with $n^{2}\alpha$ and $Q$ with
$Q  + [\ker n]$, we see that we have 
$P + Q + \Qb= [\ker\lambda] + [\ker\alpha]$. Since $\alpha$ is fixed by the 
Rosati involution associated to $\lambda$ and is totally positive, the composite
map $\nu = \lambda\alpha$ is also a polarization of~$X$, and we have 
$P + Q + \Qb= [\ker\nu]$.

Let $G = \ker\nu$ and let $e\colon G\times G\ra \BGm$ be the non-degenerate
     \marginpar{\color{Red}\small\raggedright
                Here is the error in the proof. The cited proposition
                assumes the base field is finite, but Theorem~3.2
                allows for arbitrary base fields. The highlighted
                statement is false, for example, over $\BQ$.}
alternating pairing on $G$ whose existence is shown in~\S23 
of~\cite{brefMumford}. Let $H$ be a simple element of $\Ker_\CC$  that occurs 
in~$Q$. 
     {\color{Red}Proposition~5.2 of~\cite{brefHoweAV} shows that 
     there is an embedding of $H$ into $G$ such that the pairing
     $e$ restricted to $H\times H$ is the trivial pairing.}
Let $\varphi$ be the natural isogeny from $X$ to $Y = X/H$. Then the Corollary
to Theorem~2 (p.~231) of~\S23 of~\cite{brefMumford} shows that there is a 
polarization $\nu'$ of $Y$ such that $\nu=\hat{\varphi}\nu'\varphi$. In 
$G(\Ker_{\CC})$ this gives us the equality 
$[\ker\nu] = [H] + [\ker\nu'] + \bar{[H]}.$  If we replace $Q$ by $Q-[H]$ and 
$\nu$ by~$\nu'$, we will again have the equality $P + Q + \Qb = [\ker\nu]$, but
we will have decreased the number of simple group schemes that occur in~$Q$. By 
applying this argument repeatedly, we can finally obtain the equality 
$P=[\ker\nu]$ for a polarization $\nu$ of a variety in~$\CC$. This shows that
$P$ is attainable.
\end{proof}

\goodbreak

\subsection{Theorem~\ref{TTmain} revisited}
\label{SSrevisit}
In this section we show how our proof of\break
Theorem \ref{TTmain} can be understood in terms of Theorem~\ref{TTattainable}.

Let $k$ be a $p$-admissible field, let $E$ be a $p$-isolated elliptic curve 
over~$k$, let $\ell$ be a degree-$p$ Galois extension of~$k$, let $X$ be the 
reduced restriction of scalars of $E$ from $\ell$ to~$k$, and let $\CC$ be the 
isogeny class of~$X$. Then $A = \End^{0}\CC$ is the cyclotomic field 
$K = \BQ(\zeta_{p})$, so $R_{0}(A)$ is $K^{*}$ and $R_{1}(A)$ and $R_{2}(A)$ are
both equal to the multiplicative group of totally positive elements of the
maximal real subfield $K^{+}$ of~$K$.

Note that the $E[p]$-rank defines a homomorphism $Z(\CC)/B(\CC)\ra \BZ/2\BZ$. 
Lemma~\ref{LLdegree} shows that the degree of an element of $R_{1}(A)$ is equal
to the fourth power of its norm from $K^{+}$ to~$\BQ$. Since 
Lemma~\ref{LLgaloismodule} shows that $E[p]$ is the only simple $p$-torsion 
group scheme that occurs in $\Ker_{\CC}$, and since the rank of $E[p]$ 
is~$p^{2}$, the $E[p]$-rank of every element of $\Phi(R_{1}(A))$ is even. Thus,
the $E[p]$-rank gives us a homomorphism from $\CBone(\CC)$ to $\BZ/2\BZ$. 
Lemma~\ref{LLonepolarization} shows that the image of $I_{\CC}$ under this 
homomorphism is nonzero, so $I_{\CC}$ itself is nonzero. Then 
Theorem~\ref{TTattainable} shows that the trivial group scheme is not attainable
in~$\CC$, so $\CC$ contains no principally polarized varieties.

We leave it to the reader to use the methods of this paper to prove the
following generalization of Theorem~\ref{TTmain}:

\begin{theorem}
\label{TTgeneral}
   \marginpar{\color{Red}\small\raggedright 
              A hypothesis is missing here: We should also assume
              that $p$ is not equal to the characteristic of~$k$.
              
              \quad
              
              The proof of this theorem (with the added hypothesis)
              does not use  Theorem~\ref{TTattainable} --- it is based
              on the ideas used in \S2 --- and still carries through.}
Let $\ell/k$ be a Galois extension of odd prime degree~$p$, and let $Y$ and $Z$
be abelian varieties over $k$ such that 
\begin{enumerate}
    \item $Y[p](\kb)$ is a simple $\Gal(\kb/k)$-module\textup{;}
    \item the simple module $Y[p](\kb)$ does not occur 
          in $Z[p](\kb)$\textup{;} and
    \item $\End Y_{\ell} = \BZ$.
\end{enumerate}
Let $X$ be the kernel of the trace map from $\Res_{\ell/k} Y$ to~$Y$. Then $X$ 
is simple, and every polarization of every abelian variety isogenous to 
$X\times Z$ has degree divisible by~$p^{2\dim Y}$.
\end{theorem}

\end{document}